\documentclass[10pt]{amsart}
\input{epsf}
\usepackage{amssymb,latexsym}
\usepackage{color,pstricks}
\topskip  \numberwithin{equation}{section}
\newtheorem{theorem}{Theorem}[section]

\newtheorem{corollary}[theorem]{Corollary}

\newtheorem{example}[theorem]{Example}

\begin{document}

\pagenumbering{arabic}
\pagestyle{headings}
\def\sof{\hfill\rule{2mm}{2mm}}
\def\llim{\lim_{n\rightarrow\infty}}

\title{Avoiding substrings in compositions}
\maketitle

\begin{center}
{\bf Silvia Heubach}\\
{\it Dept. of Mathematics, California State University Los
Angeles,\\
Los Angeles, CA 90032, USA}\\
{\tt sheubac@calstatela.edu}\\
\vskip 10pt
{\bf Sergey Kitaev\footnote{The work presented here was supported by grant no. 090038011 from the Icelandic Research Fund.}}\\
{\it The Mathematics Institute, School of Computer Science,
Reykjavik University, \\ 103 Reykjavik, Iceland}\\
{\tt sergey@ru.is}\\

\end{center}

\def\P{POP}
\def\A{\mathcal{A}}
\def\SS{\frak S}
\def\NN{\mathbb N}
\def\Ps{POPs}
\def\mn{\mbox{-}}
\def\newop#1{\expandafter\def\csname #1\endcsname{\mathop{\rm #1}\nolimits}}
\newop{MND}
\section*{Abstract}
A classical result by Guibas and Odlyzko obtained in 1981 gives the
generating function for the number of strings that avoid a given set
of substrings with the property that no substring is contained in
any of the others. In this paper, we give an analogue of this result
for the enumeration of compositions that  avoid a given set of
prohibited substrings, subject to the compositions' length (number
of parts) and weight. We also give examples of families of strings to be avoided that allow for an explicit formula for the generating function.  Our results extend recent results by Myers
on avoidance of strings in  compositions subject to weight,
but not length.

\noindent{\bf Keywords}: Compositions, strings, avoidance,
generating functions, (auto)correlation

\noindent{\bf 2000 Mathematics Subject Classification}: 05A05, 05A15
\thispagestyle{empty}
\section{Introduction}\label{Introduction}

\begin{figure}[htp]
\begin{center}
\begin{pspicture}(-3.5,-1)(6.6,1.65)
\put(-1.5,.15){$X_2$}
\put(-1.5,.9){$X_1$}
\psline(0,0)(3,0)(3,0.5)(0,.5)(0,0)
\pspolygon[fillstyle=solid,fillcolor=black](0,0)(1.25,0)(1.25,0.5)(0,.5)(0,0)
\put(1.55,0.15){overlap}
\psline(1.25,0.75)(4,.75)(4,1.25)(1.25,1.25)(1.25,0.75)
\put(1.55,0.9){overlap}\put(3.2,0.9){tail}
\psline[linestyle=dotted](1.25,-.2)(1.25,1.4)
\psline[linestyle=dotted](3,-.2)(3,1.4)
\end{pspicture}
\vspace{-0.35in}
\caption{Comparing strings  $X_1$ and $X_2$.}\label{Fig-tail}
\end{center}
\end{figure}
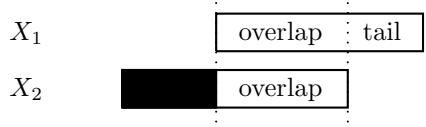

In 1981, Guibas and Odlyzko~\cite{GuibasOdlyzko} obtained
the generating function for the number of strings avoiding a given
set of prohibited substrings and then applied this result to
non-transitive games. (A string $s=s_1s_2\cdots s_m$  {\em contains a substring} $b_1b_2\cdots b_k$ of length $k$ if there is an index $i$ such that $s_is_{i+1}\cdots s_{i+k-1}=b_1b_2\cdots b_k$. Otherwise, we say that $s$ {\em avoids} the substring $b_1b_2\cdots b_k$.) A detailed derivation of this generating function  and related results in the binary case  was later given by Winterfjord in his Masters thesis~\cite{Bjorn}. The basic idea in the derivation of the generating function is the notion of the correlation between two strings and being able to enumerate the strings avoiding the set of substrings in two different ways.  Let $X_1=a_0a_1\ldots a_{m-1}$ and $X_2=b_0b_1\ldots b_{\ell-1}$ be
two strings of lengths $m$ and $\ell$, respectively, over the alphabet $[n]=\{1,2,\ldots, n\}$. The {\em correlation} $c_{12}=c_0c_1\ldots c_{m-1}$
is the binary string defined as follows: \begin{itemize}
\item[$m\leq \ell$:] For  $0\leq j\leq m-1$, $c_j=1$ if
$a_i=b_{\ell-m+i+j}$ for $i=0,1,\ldots,m-j-1$, and $c_j=0$
otherwise;
\item[$m>\ell$:] For $0\leq j\leq m-\ell$, $c_j=1$ if
$b_i=a_{m-\ell+i-j}$ for $i=0,1,\ldots,\ell-1$, and $c_j=0$
otherwise; for  $m-\ell+1\leq j\leq m-1$, $c_j=1$ if
$a_i=b_{\ell-m+i+j}$ for $i=0,1,\ldots,m-j-1$ and $c_j=0$
otherwise.
\end{itemize}
In plain English, this means that $c_j$ is equal to $1$ if and only if the coefficients in the overlap of the string $X_1$ and the string $X_2$, shifted (or offset) by $j$ positions to the left, agree, as illustrated in Figure~\ref{Fig-tail}.

For example, if
$X_1=110$ and $X_2=1011$, then $c_{12}=011$ and $c_{21}=0010$, as
depicted below:
$$\begin{array}{c|cccccccc}
 \text{offset }j & &   &   & 1 & 1 & 0 & \ &c_j   \\
\hline
0& &   & 1 & 0 & 1 & 1 &   & 0 \\
 1& & 1 & 0 & 1 & 1 &   &   & 1 \\
2&1 & 0 & 1 & 1 &   &   &   & 1 \\
& & & & & & & & \\
\end{array}
\begin{array}{cc}
\ & \ \\
\end{array}
\begin{array}{c|cccccccc}
 \text{offset }j &  &   & 1 & 0 & 1 & 1 & \ & c_j  \\
\hline
0&  &   &   & 1 & 1 & 0 &   & 0 \\
1&  &   & 1 & 1 & 0 &   &   & 0 \\
2&  & 1 & 1 & 0 &   &   &   & 1 \\
3&1 & 1 & 0 &   &   &   &   & 0 \\
\end{array}$$
In general $c_{12}\neq c_{21}$ and, unless the strings are of the same lengths, the correlations will have different lengths. The {\em autocorrelation} of a string or word $X_1$ is just
$c_{11}$, the correlation of $X_1$ with itself. For instance, if
$X_1=1011$ then $c_{11}=1001$. It is convenient to associate
 a  {\em correlation polynomial}
$c_{12}(x)=c_0+c_1q+\cdots +c_{k-1}q^{k-1}$ with the correlation $c_{12}=c_0c_1\ldots c_{k-1}$. This correlation polynomial is the generating function for the number of letters in the {\em tail},  the portion that is to the right of the overlap in the substring $X_1$, as illustrated in Figure~\ref{Fig-tail}.

We now state the general result given by Guibas and Odlyzko~\cite{GuibasOdlyzko} in the form given (for the special case of binary strings) in Winterfjord \cite[Th. 24]{Bjorn}.

\begin{theorem}\label{stringgf}
The generating function for the number of strings or words of length $n$ over a given alphabet
 that avoid the substrings $S_1,\ldots, S_k$ of lengths
$\ell_1, \ldots,\ell_k$ respectively, none included in any other, is
given by
\begin{equation}\label{main} S(q)=
\frac{\begin{array}{|ccc|} -c_{11}(q) & \cdots & -c_{1n}(q) \\
\vdots & \ddots & \vdots \\ -c_{n1}(q) & \cdots & -c_{nn}(q)
\end{array}}{\ \ \begin{array}{|cccc|} (1-nq) & 1 & \cdots & 1 \\ q^{\ell_1} & -c_{11}(q) & \cdots & -c_{1n}(q) \\
 \vdots & \vdots & \ddots & \vdots \\ q^{\ell_k} & -c_{k1}(q) & \cdots & -c_{kk}(q) \end{array}\ \ }\ , \end{equation}
where $c_{ij}(q)$ is the correlation polynomial for the substrings $S_i$ and $S_j$.\end{theorem}

 Unfortunately, the approach by Guibas and Odlyzko is not applicable to
permutations and subpermutations, or when {\em patterns} (as opposed
to strings) are to be avoided. However, the approach generalizes to
{\em compositions} avoiding a set of prohibited substings, and we
will derive a formula for the most general case that is an analogue
of the formula by Guibas and Odlyzko\footnote{As the matter of
fact, a recent paper by Myers~\cite{Myers} considers a very similar
problem. However, we are able to control
both length and weight in compositions, as
opposed to just weight, while Myers' result is more general with respect to the alphabet considered.}. This generalization to
compositions follows the current interest in compositions which have
been studied from different perspectives in the literature, mostly
from the view point of pattern avoidance (see~\cite{HM} and
references therein). Our results add a facet to this research.

Let $\NN$ be the set of natural numbers. A {\em composition} $\sigma=\sigma_1\cdots\sigma_m$ of $n\in\NN$ is an ordered collection (or string) of one or more
positive integers whose sum, also called the composition's {\em
weight} $w(\sigma)$, is $n$. The number of {\em summands} or {\em
letters}, namely $m$, is called the number of {\em parts} of the
composition and is denoted by $\ell(\sigma)$.
 The main result of this paper is the derivation of  the generating function
$$G(x,q)=G(S_1,\ldots,S_k;x,q)=\sum_{\sigma}x^{w(\sigma)}q^{\ell(\sigma)}$$
where the sum is taken over all compositions with parts in $\NN$  simultaneously avoiding
the prohibited substrings $S_i$, $i=1,\ldots,k$, where none of the substrings is included in any other. We state and prove this result in Section~\ref{compositions} and then give applications of our result  for families
of prohibited substrings in Section~\ref{applications}.

\section{Main result }\label{compositions}

In order to generalize Theorem~\ref{stringgf} to compositions, we need to adapt the correlation polynomial to also keep track of the the weight in addition to the length of the tail.  We therefore define the correlation polynomial for a correlation $c_{ij}=c_0c_1\ldots c_{m-1}$ between $S_i=a_0a_1\ldots
a_{m-1}$ and $S_j$ as
$$c_{ij}(x,q)=c_0+c_1x^{w(a_{m-1})}q+c_2x^{w(a_{m-2}a_{m-1})}q^2+\cdots
+c_{m-1}x^{w(a_{2}a_{3}\ldots a_{m-1})}q^{m-1}.$$  For example, for $X_1=110$ and
$X_2=1011$ considered in Section~\ref{Introduction},
$c_{12}(x,q)=x+x^2q$, $c_{21}(x,q)=(xq)^2$, $c_{11}(x,q)=1$, and
$c_{22}(x,q)=1+x^3q^2$. Note that since we are considering compositions, all parts are positive and therefore each term but the first one of a
correlation polynomial  is divisible by $x q$
(the first term is either 0 or 1). We are now ready to state the main result.

\begin{theorem}\label{mainthm}
The generating function for the number of compositions of weight $n$
and  length $m$ with parts in $\NN$ that avoid the substrings
$S_1,\ldots, S_k$ of lengths $\ell(S_1), \ldots,\ell(S_k)$
respectively, none included in any other, is given by

\begin{equation}\label{main-res} G(x,q)=
\frac{(1-x)\cdot\begin{array}{|ccc|} -c_{11}(x,q) & \cdots & -c_{1n}(x,q) \\
\vdots & \ddots & \vdots \\ -c_{n1}(x,q) & \cdots & -c_{nn}(x,q)
\end{array}}{\ \ \begin{array}{|cccc|} 1-x(1+q) & 1-x & \cdots & 1-x \\ x^{w(S_1)}q^{\ell(S_1)} & -c_{11}(x,q) & \cdots & -c_{1n}(x,q) \\
 \vdots & \vdots & \ddots & \vdots \\ x^{w(S_k)}q^{\ell(S_k)} & -c_{n1}(x,q) & \cdots & -c_{nn}(x,q) \end{array}\ \ }\  \end{equation}
where  $c_{ij}(x,q)$ are the correlation polynomials defined above.
\end{theorem}

\begin{proof}
In finding $G(x,q)$ we adapt the arguments in~\cite{GuibasOdlyzko,Bjorn} to compositions. Let $A$ denote the set of all compositions avoiding the prohibited
substrings and let
$B_i$, for $i=1,\ldots,k$, be the set of all compositions
ending with $S_i$ but having no other occurrence of any of the prohibited
substrings. A composition in $B_i$ is said to {\em quasi-avoid} $S_i$. We denote the generating function corresponding to $B_i$ by $B_i(x,q)$ and note that $G(x,q)$ is the generating function of the set $A$. Furthermore, the sets $A$ and $B_i$ are all pairwise disjoint as none of the substrings is included in any of the others.

We now derive recurrences for certain sets of compositions. Note that we can create compositions of weight $n+1$ recursively from those of weight $n\ge1$ by either increasing the last part by $1$ or by appending a part $1$ at the right end of the composition. For a set of compositions $M$, let $M^{+1}$ denote the set obtained
from $M$ by increasing the rightmost part of {\em each} non-empty composition
by $1$, and let $M\times \{1\}$ denote the set obtained from $M$ by
adjoining the new rightmost part $1$ to {\em each} composition in $M$.
With this notation, we can express the set of compositions that either avoid or quasi-avoid the substrings  as follows:
\begin{equation}\label{set1} A\cup
B_1\cup\cdots\cup B_k=\{\epsilon\}\cup (A\cup B_1\cup \cdots\cup
B_k-\{\epsilon\})^{+1}\cup (A\times\{1\}),\end{equation} where $\epsilon$ is the empty composition. The expression on the right hand side follows as increasing the last part of a composition that avoids all substrings can create an occurrence of a substring, but  only at the very end of the composition, and likewise when adding a new part. On the other hand, a composition that  quasi-avoids a string is transformed either into a composition that avoids the substrings or quasi-avoids a different substring when increasing the last part. However, when appending the part $1$ to a composition that quasi-avoids $S_i$ we create a composition that contains $S_i$, so that operation is not allowed for the sets $B_i$. Increasing the last part results in an increase in the weight of the composition by $1$ but no increase in the number of parts, while appending the part $1$ increases both the weight and the length of the composition. Thus  \eqref{set1} can be expressed in terms of generating functions as \begin{equation}\label{first}
(1-x -x q)G(x,q)+(1-x)(B_1(x,q)+\cdots+B_k(x,q))=1-x,
\end{equation}
where we have used that the generating function of the union of disjoint sets is the sum of the respective generating functions, and the generating function of a Cartesian product is the product of the respective generating functions.

We now create an alternative connection between the sets $A$ and $B_i$. Let $R_i$ denote the set of compositions that consist of a composition
from $A$ followed by the prohibited string $S_i$, where
$i=1,\ldots,k$. Note that $R_i$ and $R_j$ are disjoint for $i\neq j$
as none of $S_i$'s is included in any other. Furthermore, the set $R_i$ is not
identical to the set $B_i$ as it is possible that a prohibited
string will occur {\em inside} a string in $R_i$, not only
at the end. For a composition (or string) $X$ from $B_j$, we call a string $Y$ with
$\ell(Y)\leq \ell(S_i)-1$ a {\em possible $ij$-tail} if $XY$ ends
with the substring $S_i$. This nomenclature is readily understood when comparing Figure~\ref{Fig-ijtail} to Figure~\ref{Fig-tail}, as $Y$ is the tail in the comparison of $S_i$ with $S_j$.

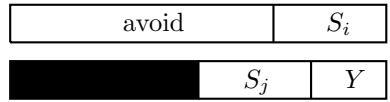
\begin{figure}[htp]
\begin{center}
\begin{pspicture}(-1,0)(4,1.5)
\psline(1.5,0)(4,0)(4,0.5)(1.5,0.5)(1.5,0)
\psline(3,0)(3,0.5)
\pspolygon[fillstyle=solid,fillcolor=black](-1,0)(1.5,0)(1.5,0.5)(-1,.5)(-1,0)
\put(2.1,.15){$S_j$}\put(3.45,0.15){$Y$}
\psline(-1,.75)(4,.75)(4,1.25)(-1,1.25)(-1,0.75)
\psline(2.5,.75)(2.5,1.25)
\put(0.5,.9){avoid}\put(3.2,.9){$S_i$}
\end{pspicture}
\caption{The $ij$-tail $Y$.}\label{Fig-ijtail}
\end{center}
\end{figure}
\vspace{-0.2in}

With this definition, we obtain the following equality of sets:
\begin{equation}\label{set2} A\times S_i=\cup_{1\leq j\leq
k}B_j\times \{\mbox{possible $ij$-tail}\},\end{equation} which in
terms of generating functions gives the following equation for each
$i=1,\ldots,k$:
\begin{equation}\label{second}
G(x,q)x^{w(S_i)}q^{\ell(S_i)}-\sum_{j=1}^kc_{ij}(x,q)B_j(x,q)=0.
\end{equation}
Indeed, a proof of~(\ref{set2}) is identical to the corresponding
statement for strings that can be found
in~\cite{GuibasOdlyzko,Bjorn} (it does not matter whether we deal
with strings or compositions in this case), while for the generating
functions, the difference is that we also keep track of the weight in the
compositions using the variable $x$.

Combining~(\ref{first}) and~(\ref{second}) results in the following set of equations
$$\left(\begin{array}{cccc} 1-x(1+q) & 1-x & \cdots & 1-x \\
x^{w(S_1)}q^{\ell(S_1)} & -c_{11}(x,q) & \cdots & -c_{1n}(x,q) \\
 \vdots & \vdots & \ddots & \vdots \\
 x^{w(S_k)}q^{\ell(S_k)} & -c_{n1}(x,q) & \cdots & -c_{nn}(x,q) \end{array}\right)
 \left(\begin{array}{c}G(x,q) \\ B_1(x,q)\\ \vdots \\
 B_k(x,q)\end{array}\right)=\left(\begin{array}{c}1-x \\ 0\\ \vdots \\
 0\end{array}\right).$$
Using Cramer's rule to solve for $G(x,q)$ gives formula~(\ref{main-res}).
\end{proof}

\section{Applications of Theorem~\ref{mainthm}}\label{applications}

Even though Theorem~(\ref{mainthm}) provides an explicit solution
to the enumerative problem, it involves the evaluation of determinants
which may not be a simple thing to do. While one can easily find explicit formulas for the
generating function that do not involve determinants when there are
just a few prohibited substrings, it is interesting to know in which
cases the determinants can be evaluated for families of
prohibited substrings. In this section, we evaluate the determinants
for a family of prohibited substrings which generalizes the
 {\em well-based sets} used in~\cite{Kit} to count independent
sets in certain graphs called {\em path-schemes}.

Let $1^i$ denote the string consisting of $i$ $1$'s and let $V=\cup_{1\leq i\leq k}\{21^{a_i-1}2\}$ with $1\leq a_1<a_2<\cdots
<a_k$ be the  set of substrings to be avoided. Note that none of the
substrings in $V$ is included in any other. Thus
we can apply formula~(\ref{main-res}) to find the generating
function for the number of compositions avoiding all the substrings
in $V$ simultaneously.

\begin{corollary} The generating function $V(x,q)$ for the number of compositions of weight $n$ and  length $m$ with parts in $\NN$ that avoid the family of substrings $V$ defined above is
given by
\begin{equation}\label{main-appl}
V(x,q)=\frac{(1-x)(1+x\sum_{i=1}^{k}(xq)^{a_i})}{(1-x(1+q)+(1-x)x^2q)(1+x\sum_{i=1}^{k}(xq)^{a_i})-(1-x)x^2q}.
\end{equation}
\end{corollary}

\begin{proof} It is easy to see
that the correlation polynomial for the two strings $21^{a_i-1}2$ and $21^{a_j-1}2$
is $c_{ij}(x,q)=\delta_{ij}+x(xq)^{a_i}$, where $\delta_{ij}$ is the
Kronecker delta. Also,
$x^{w(21^{a_i-1}2\})}q^{\ell(21^{a_i-1}2)}=x^{a_i+3}q^{a_i+1}$. Therefore Theorem~\ref{mainthm} gives  that
$$V(x,q)= \frac{(1-x)\cdot\begin{array}{|cccc|} -1-x(xq)^{a_1} &
-x(xq)^{a_1} &  \cdots
&  -x(xq)^{a_1} \\
-x(xq)^{a_2} &  -1-x(xq)^{a_2} &  \cdots
&  -x(xq)^{a_2} \\
\vdots & \vdots & \ddots & \vdots \\
-x(xq)^{a_k} &  -x(xq)^{a_k} &  \cdots &  -1-x(xq)^{a_k}
\end{array}}{\ \ \begin{array}{|ccccc|}
1-x(1+q) & 1-x & 1-x & \cdots & 1-x \\
x^{a_1+3}q^{a_1+1} &  -1-x(xq)^{a_1} &  -x(xq)^{a_1} &  \cdots
&  -x(xq)^{a_1} \\
x^{a_2+3}q^{a_2+1} & -x(xq)^{a_2} &  -1-x(xq)^{a_2} &  \cdots
&  -x(xq)^{a_2} \\
 \vdots & \vdots & \vdots & \ddots & \vdots \\
x^{a_k+3}q^{a_k+1} & -x(xq)^{a_k} &  -x(xq)^{a_k} &  \cdots &
-1-x(xq)^{a_k}  \end{array}\ \ }.$$

To compute the determinant in the numerator, replace row 1 by the
sum of all rows and then factor out the common factor
$(-1-x\sum_{i=1}^{k}(xq)^{a_i})$. Next subtract column 1 from
columns $2,3,\ldots,k$ to obtain
$$ -(1+x\sum_{i=1}^{k}(xq)^{a_i}) \cdot
\begin{array}{|ccccc|} 1 & 0 & 0 & \cdots & 0 \\
-x(xq)^{a_2} & -1 & 0 & \cdots & 0 \\ -x(xq)^{a_3} & 0 & -1 & \cdots & 0 \\
\vdots & \vdots & \vdots & \ddots & \vdots \\ x(xq)^{a_k} & 0 & 0 &
\cdots & -1 \end{array}=(-1)^k\cdot (1+x\sum_{i=1}^{k}(xq)^{a_i}).$$

To compute the determinant in the denominator, replace column 1 by
the sum of column 1 and  $x^2q\cdot$(column ($k+1$)) and for $i=2,3,\ldots,k$, replace
column $i$ by the difference of column
$i$ and (column ($k+1$)) to yield
$$\begin{array}{|ccccc|}
1-x(1+q)+(1-x)x^2q & 0 & \cdots & 0 & 1-x\\
0 &  -1 & \cdots  &  0 &  -x(xq)^{a_1} \\
0 & 0 &  \cdots &  0 &  -x(xq)^{a_2} \\
 \vdots & \vdots & \ddots & \vdots & \vdots \\
0 & 0 &  \cdots &  -1 & -x(xq)^{a_{k-1}} \\
-x^2q & 1 &  \cdots &  1 & -1-x(xq)^{a_{k}} \\
\end{array}\ .$$

To obtain an upper triangular matrix we replace the last row in this determinant by
$$\frac{x^2q\mbox{(row 1)}}{1-x(1+q)+(1-x)x^2q}+\mbox{(row 2)}+\mbox{(row 3)}+\cdots+
\mbox{(row }(k+1))$$
which yields that the determinant of the denominator is equal to
$$(-1)^k\Big[(1-x)x^2q-(1-x(1+q)+(1-x)x^2q)(1+x\sum_{i=1}^{k}(xq)^{a_i})\Big], $$
completing the proof.
\end{proof}

Further simplifications of $V(x,q)$ are possible whenever
$\sum_{i=1}^{k}(xq)^{a_i}$ can be simplified. We provide three examples
here.

\begin{example}
The set of prohibited substrings $\{22, 212, \ldots, 2i^{k-1}2\}$ corresponds to $\{a_1,a_2,\ldots,a_k\}=\{1,2,\ldots,k\}$. In this
case,~(\ref{main-appl}) reduces to
$$V_k(x,q)=\frac{(1-x)(1-x q+x^2 q (1-(x q)^k))}{(1-x (1+q)+(1-x)x^2q)(1-x q+x^2 q (1-(x q)^k)) -(1-x)(1-x q)x^2q }.$$

The initial values of $V_2(x,q)$ \rm{(}avoiding $22$ and $212$\rm{)} are as follows:
\begin{eqnarray*}V_2(x,q) &=& 1+q x+\left(q+q^2\right) x^2+\left(q+2 q^2+q^3\right) x^3+\left(q+2 q^2+3 q^3+q^4\right) x^4+\\
&& \left(q+4 q^2+3 q^3+4 q^4+q^5\right) x^5+\left(q+5 q^2+9 q^3+5
q^4+5 q^5+q^6\right) x^6 +\cdots
\end{eqnarray*}
\end{example}

\begin{example} The set of prohibited substrings that have an even number of $1$'s, $\{22, 2112, \ldots, 2i^{2k}2\}$ is represented by the set
$\{a_1,a_2,\ldots\}=\{1,3,5,\ldots,2k+1\}$. In this
case,~(\ref{main-appl}) is simplified as follows:
$$V_o(x,q)=\frac{(1-x) \left(1-(xq)^2+x^2q \left(1-(xq)^{2
k+1}\right)\right)}{\left(1-(1+q) x+(1-x) x^2q\right)
\left(1-(xq)^2+xq^2 \left(1-(xq)^{2 k+1}\right)\right)-(1-x) x^2q
\left(1- (xq)^2\right)}.$$ The initial values of $V_o(x,q)$ for
$k=2$  \rm{(}avoiding $\{22, 2112, 211112\}$\rm{)} are as follows: \begin{eqnarray*}V_o(x,q) &=& 1 + xq + (q +
q^2) x^2 + (q + 2 q^2 + q^3) x^3 + (q + 2 q^2 + 3 q^3 +
     q^4) x^4 + \\ && (q + 4 q^2 + 4 q^3 + 4 q^4 + q^5) x^5 + (q + 5 q^2 +
    9 q^3 + 6 q^4 + 5 q^5 + q^6) x^6 + \\ && (q + 6 q^2 + 13 q^3 + 16 q^4 +
    9 q^5 + 6 q^6 + q^7) x^7 + \\
    && (q + 7 q^2 + 19 q^3 + 28 q^4 + 26 q^5 + 12 q^6 + 7 q^7 + q^8) x^8 + \cdots.
    \end{eqnarray*}
    \end{example}

\begin{example} The set of prohibited substrings that have an odd number of $1$'s, $\{212, 21112, \ldots, 2i^{2k-1}2\}$ is represented by the set $\{a_1,a_2,\ldots\}=\{2,4,6,\ldots,2k\}$. In this
case,~(\ref{main-appl}) is simplified as follows:
$$V_e(x,q)=\frac{(1-x) \left(1-(xq)^2+x^3q^2 \left(1-(xq)^{2 k}\right)\right)}
{\left(1-(1+q) x+(1-x) x^2q\right) \left(1-(xq)^2+x^3q^2
\left(1-(xq)^{2 k}\right)\right)-(1-x) x^2q
\left(1-(xq)^2\right)}.$$ The initial values of $V_e(x,q)$ for
$k=2$  \rm{(}avoiding $\{212, 21112\}$\rm{)} are as follows:
\begin{eqnarray*}V_e(x,q) &=& 1+xq+\left(q+q^2\right) x^2+\left(q+2 q^2+q^3\right) x^3+\left(q+3 q^2+3 q^3+q^4\right) x^4+\\
&& \left(q+4 q^2+5 q^3+4 q^4+q^5\right) x^5+\left(q+5 q^2+10 q^3+8
q^4+5 q^5+q^6\right) x^6+\\ && \left(q+6 q^2+15 q^3+18 q^4+11 q^5+6
q^6+q^7\right) x^7+\\ && \left(q+7 q^2+21 q^3+33 q^4+30 q^5+15
q^6+7q^7+q^8\right) x^8 + \cdots.
\\  \end{eqnarray*}
\end{example}

Clearly, other families of substrings can be created that allow for similar simplification of the generating function.


\def\dtcs{{\em Discrete Mathematics and Theoretical Computer Science\,}}
\def\pias{{\em Procedure Indian Academic Science\,}}
\def\eujc{{\em European J. Combin.\,}}
\def\aam{{\em Advances in Applied Mathematics\,}}
\def\dm{{\em Discrete Math.\,}}
\def\ejc{{\em Electron. J. Combin.\,}}
\def\jcta{{\em Journal Comb. Theory Series A\,}}
\def\slc{{\em S\'eminaire Lotharingien de Combinatoire\,}}
\def\ars{{Ars Combinatorica\,}}

\def\vv{{Volume\,}}

\end{document}